# Specification-Oriented Automatic Design of Topologically Agnostic Antenna Structure

Adrian Bekasiewicz[1][0000-0003-0244-541X], Mariusz Dzwonkowski[1,2][0000-0003-3580-7448], Tom Dhaene[3][0000-0003-2899-4636], and Ivo Couckuyt[3][0000-0002-9524-4205]

[1] Faculty of Electronics, Telecommunications and Informatics, Gdansk University of Technology, Narutowicza 11/12, 80-233 Gdansk, Poland
[2] Department of Radiology Informatics and Statistics, Faculty of Health Sciences, Medical University of Gdansk, Tuwima 15, 80-210 Gdansk, Poland.
[3] Department of Information Technology (INTEC), IDLab, Ghent University-imec, iGent, Technologiepark-Zwijnaarde 126, 9052 Ghent, Belgium
`bekasiewicz@ru.is`

**Abstract.** Design of antennas for modern applications is a challenging task that combines cognition-driven development of topology intertwined with tuning of its parameters using rigorous numerical optimization. However, the process can be streamlined by neglecting the engineering insight in favor of automatic determination of structure geometry. In this work, a specification-oriented design of topologically agnostic antenna is considered. The radiator is developed using a bi-stage algorithm that involves min-max classification of randomly-generated topologies followed by local tuning of the promising designs using a trust-region optimization applied to a feature-based representation of the structure frequency response. The automatically generated antenna is characterized by –10 dB bandwidth of over 600 MHz w.r.t. the center frequency of 6.5 GHz and a dual-lobe radiation pattern. The obtained performance figures make the radiator of use for in-door positioning applications. The design method has been favorably compared against the frequency-based trust-region optimization.

## 1 Introduction

Design of modern antennas is an inherently cognitive process. It involves experience-driven development of topology followed by its tuning so as to fulfill the desired performance requirements [1], [2], [3]. When combined with robust optimization, this engineer-in-a-loop approach proved to be useful for the design of new (often unconventional) topologies [4], [5]. Although restricting the antenna development to a specific shape is considered pivotal for ensuring feasibility of the design process, it is also subject to engineering bias and hence limits the potential in terms of achieving solutions characterized by unique (perhaps not expected) performance characteristics [4]. These might include, for instance broadband/multiband behavior, but also im-proved radiation capabilities and/or small dimensions [3]-[5]. From this perspective, a streamlined design that rely only on numerical optimization methods represents an interesting alternative to the standard antenna development techniques [6].

Automatic antenna generation governed by numerical optimization is a challenging task. From the perspective of geometry, the design can be represented as a set of points (interconnected using, e.g., line-sections, or splines), or in the form of a binary matrix that defines the configuration of primitives (e.g., rectangles) constituting the antenna [6]-[11]. Point-based approaches are capable of supporting geometrically complex topologies. Owing to continuous and free-form nature of coordinates, "evolution" of geometry can be governed by standard numerical optimization algorithms [7], [9]. However, a large number of points and constraints on their distribution are required to ensure flexibility in terms of attainable topology and its consistency (here understood as lack of self-intersections of the coordinate-based curves) [7], [8]. Matrix-based methods naturally ensure consistency by representing topologies as compositions of partially overlapping primitives [6], [10]. Their main bottleneck is a large number of variables required even for relatively simple topologies. Besides activation/deactivation of primitives (based on contents of the matrix), their dimensions also need to be scaled which reveals a mixed-integer nature of the problem [10]. Finally, high dimensionality and the need to evaluate performance based on the expensive electromagnetic (EM) simulations make the discussed universal representations impractical for optimization when conventional algorithms are considered, as they require hundreds or even thousands of design evaluations to converge [4], [9], [10].

The problem related to unacceptable numerical cost can be mitigated using trust-region (TR) methods. TR is a class of techniques that enable iterative approximation of the desired solutions based on the optimization of data-efficient models [12], [13]. Regardless of proved usefulness, local character of the TR-based methods limit their usefulness for improving already promising initial designs rather than exploration of the feasible space. Furthermore, due to complexity and multi-modal character of the problem, successful development of antennas using the outlined routines is a subject to availability of appropriate response-processing mechanisms.

In this work, a framework for specification-oriented development of topologically agnostic antenna has been considered. The method involves classification of automatically generated topologies followed by TR-based tuning of the selected design in a setup where the frequency responses are represented in the form of carefully selected feature points so as to enable the algorithm convergence. The method has been used to design a 53-dimensional planar antenna represented using a set of points interconnected using line sections. The structure has been optimized so as to ensure reflection below –10 dB over a frequency range from 6.2 GHz to 6.8 GHz which corresponds to 5th channel of the ultra-wideband (UWB) spectrum. The antenna is also characterized by a dual-lobe radiation pattern and a relatively high gain of up to 7 dB. Owing to planar topology, the proposed radiator might be used as a component of the in-door localization systems, or medical-imaging devices (e.g., for breast cancer detection). Comparison of the considered design approach against conventional, TR-based optimization over the frequency domain has also been considered.

## 2    Topologically Agnostic Antenna

Automatic specification-driven antenna development is subject to availability of a model that supports free-form adjustment of geometries. Consider a topology-agnostic planar antenna shown in Fig. 1. The structure is implemented on a substrate with permittivity and thickness of $\varepsilon_r = 2.55$ and $h = 1.52$ mm, respectively. It comprises a patch fed through a concentric probe. The radiator is represented using the interconnected points defined in a cylindrical coordinate system. To ensure feasibility of geometry each angular coordinate is set as relative to the previous one (the first is set to 0) and all angles are scaled w.r.t. their cumulative sum and normalized to the range between 0 and $2\pi$ (full circle). It should be noted that symmetry planes for the radiator are not defined so as to ensure its flexibility in terms of attainable performance specifications (even though this is achieved at the expense of increased computational cost associated with a large number of design parameters). Note that such representation supports random generation of topologies. For the probe, an additional routine is used to ensure that it is enclosed within the generated shape.

The antenna is implemented in a CST Microwave Studio and evaluated using its time-domain simulator [14]. To ensure flexibility in terms of the number of design parameters, the structure EM model is generated dynamically and implements a safeguard mechanisms that verify the consistency of geometry and handle errors [14]. On average, the EM models are discretized using 100,000 tetrahedral mesh cells and the cost of single evaluation amounts to 60 s.

The antenna geometry is represented using the following vector of design parameters $\boldsymbol{x} = [C\ \rho_f\ \varphi_f\ \boldsymbol{\rho}\ \boldsymbol{\varphi}]^T$, where $C$ is a scaling coefficient for the patch/probe coordinates, $\rho_f$ and $\varphi_f$ represent the radial and angular position of the feed w.r.t. the origin of the cylindrical system, whereas the vectors $\boldsymbol{\rho} = [\rho_1\ ...\ \rho_l]^T$ and $\boldsymbol{\varphi} = [\varphi_1\ ...\ \varphi_l]^T$ ($l = 1, ..., L$) comprise the points that represent the shape of the patch (cf. Fig. 1). The parameters $o = 5$, $r_1 = 1.27$, and $r_2 = 2.84$ remain fixed, whereas $A = B = 2(C \cdot \max(\boldsymbol{\rho}) + o)$. Note that the specific values of $r_1$ and $r_2$ ensure 50 Ohm input impedance of the feed. All dimensions except for the fixed ones and $C$ (in mm) are unit-less. The antenna topology is considered feasible within the following lower and upper design bounds: $\boldsymbol{l} = [25\ 0\ \varphi_{0.l}\ 0.1\boldsymbol{1}\ 0.01\boldsymbol{1}]^T$ and $\boldsymbol{u} = [35\ \rho_{0.f}\ \varphi_{0.h}\ 0.9\boldsymbol{1}\ 0.8\boldsymbol{1}]^T$, where $\boldsymbol{1}$ is $L$-dimensional vector of ones; $\varphi_{0.l} = \varphi^{(0)} - \pi$, $\varphi_{0.h} = \varphi^{(0)} + \pi$, and $\rho_{0.f} = \max(\boldsymbol{\rho}^{(0)})$ are the constraints determined based on the randomly generated design $\boldsymbol{x}^{(0)}$ considered for numerical optimization based on the classification (cf. Section 3). Note that the vector $\boldsymbol{x}$ comprises $D = 2L + 3$ design parameters.

## 3    Design Methodology

### 3.1    Problem Formulation

Let $\boldsymbol{R}(\boldsymbol{x}) = \boldsymbol{R}(\boldsymbol{x}, \boldsymbol{f})$ be the antenna reflection response obtained over the frequency sweep $\boldsymbol{f}$ for the vector of input parameters $\boldsymbol{x}$. The design problem can be formulated as the following non-linear minimization problem:

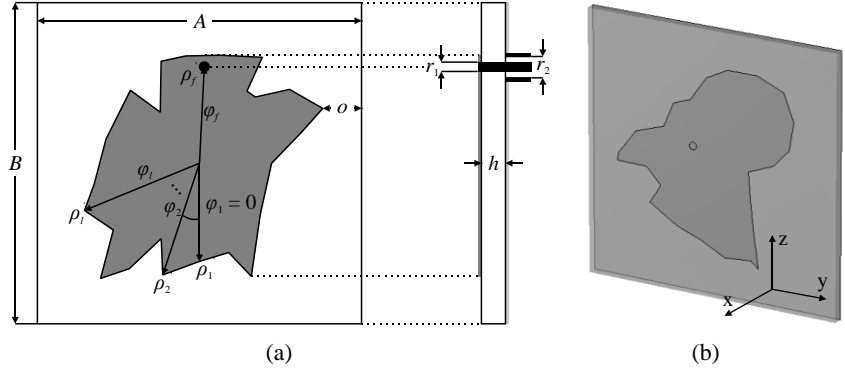

(a)           (b)

**Fig. 1.** The proposed topologically agnostic antenna: (a) geometry of the structure with highlight on the design parameters and (b) visualization of the structure optimized using the algorithm of Section 3. Note that $\rho_f' = C\rho_f$ and $\rho_l' = C\rho_l$ ($l = 1, \ldots, L$).

$$\boldsymbol{x}^* = \arg\min_{\boldsymbol{x}\in X}\left(U\left(\boldsymbol{R}(\boldsymbol{x})\right)\right) \qquad (1)$$

where $\boldsymbol{x}^*$ is the optimal design to be found within the feasible region of the search space $\boldsymbol{X}$ and $U$ is a scalar objective function. To enable unsupervised design of topology-agnostic antennas the problem (1) is solved using a bi-stage procedure where the promising candidates are first identified and then optimized using a gradient algorithm embedded within a TR framework.

The problem concerning automatic determination of the promising initial design is addressed through a random generation of the candidate solutions followed by their sorting based on a simple min-max classifier. Let $\boldsymbol{X}_r = \{\boldsymbol{x}_n\}$, $n = 1, \ldots, N$, ($\boldsymbol{X}_r \subset \boldsymbol{X}$) be an $N$-element set of random designs that represent the geometrically-flexible antenna and $\{\boldsymbol{R}_n(\boldsymbol{x}_n)\}$ be the set of their corresponding EM-based responses. Each of the obtained designs can be evaluated as:

$$E_n = \max\left(\boldsymbol{R}(\boldsymbol{x}_n)\right)_{f_l \leq f \leq f_H} \qquad (2)$$

where $f \in \boldsymbol{f}$ represent the frequency points within the bandwidth from $f_l$ to $f_h$. Upon evaluation, the responses $E_1, \ldots, E_N$ are sorted according to their values and the design $\boldsymbol{x}_k \in \boldsymbol{X}_r$ that corresponds to the solution featuring the lowest value is selected as a starting point for local TR-based optimization.

### 3.2 Feature-Based Representation of Responses

The main challenge pertinent to specification-oriented design of topology-agnostic structures is that due to intricate relations between the input parameters, the resulting responses are non-linear functions of frequency. The problem can be mitigated by shifting the design problem to a feature-based domain where the key properties of the frequency response (from the optimization perspective) are represented using a set of carefully selected points.

Let $F(x) = P(R(x))$ be the response of the antenna expressed in terms of the feature points, where $P$ is the function used for their extraction. The feature-based response is defined as $F(x) = [\omega\ S]^T$, where $\omega = [\omega_1 \ldots \omega_q]^T$ and $S = [S_1 \ldots S_q]^T$ ($q = 1, \ldots, Q$) represent the specific points related to frequencies of interest and response levels (cf. Fig. 2(a)). In other words, a pair of points $\omega_\gamma$, $S_\gamma$ ($\gamma \leq Q$) can represent specific frequency (e.g., pertinent to a local minimum of the response) and its corresponding reflection level, whereas $\omega_\kappa$, $S_\kappa$ ($\kappa \leq Q$; $\kappa \neq \gamma$) can refer to the frequency point at a specific pre-determined reflection level (e.g., at the edge of the bandwidth). Note that although a pair of points is stored, one typically requires a specific component of the pair (e.g., pertinent to resonant frequency, or level at a local maximum). Compared to the original response, its feature-based representation is a much less-nonlinear function of input parameters and hence aids design optimization (especially when performed using linear approximation models). The outlined concept is illustrated in Fig. 2. For more comprehensive discussion on the feature-based representation of frequency responses, see [15], [16].

### 3.3  Optimization Algorithm

The antenna is optimized using a gradient-based algorithm embedded in a trust-region framework. The latter generates approximations, $i = 0, \ldots$, to the final design by solving [12]:

$$x^{(i+1)} = \arg \min_{\|x - x^{(i)}\| \leq \delta} \left( U\left( G_\varepsilon^{(i)}(x) \right) \right) \quad (3)$$

The composite model $G_\varepsilon^{(i)} = [G_\omega^{(i)}\ G_S^{(i)}]^T$, where $G_\omega^{(i)} = \omega(x^{(i)}) + J_\omega(x^{(i)})(x - x^{(i)})$ and $G_S^{(i)} = S(x^{(i)}) + J_S(x^{(i)})(x - x^{(i)})$ represents the first-order surrogates generated around $x^{(i)}$ w.r.t. the frequency- and level-related features. The Jacobians are based on a large-step finite differences (FD) [15]:

$$J_\omega\left(x^{(i)}\right) = \left[ \cdots \quad \left(\omega\left(x^{(i)} + p_d^{(i)}\right) - \omega\left(x^{(i)}\right)\right)\frac{1}{p_d^{(i)}} \quad \cdots \right]^T$$
$$J_S\left(x^{(i)}\right) = \left[ \cdots \quad \left(S\left(x^{(i)} + p_d^{(i)}\right) - S\left(x^{(i)}\right)\right)\frac{1}{p_d^{(i)}} \quad \cdots \right]^T \quad (4)$$

Note that the parameter $p_d^{(i)}$ ($d = 1, \ldots, D$) denotes the FD perturbation w.r.t. $d$th dimension of the (currently best) design $x^{(i)}$, whereas the components of vector $p_d^{(i)}$ are set to zero for all dimensions except $d$th which is equal to $p_d^{(i)}$. The FD steps are selected in proportion to the design $x^{(i)}$ and updated after each successful iteration [17]. The radius $\delta$ is controlled based on the ratio $\rho = [U(F(x^{(i+1)})) - U(F(x^{(i)}))]/[U(G_\varepsilon^{(i)}(x^{(i+1)})) - U(G_\varepsilon^{(i)}(x^{(i)}))]$. The initial radius is set to $\delta = 1$. When $\rho < 0.25$ (poor prediction of $G_\varepsilon^{(i)}$) then $\delta = \delta/3$, whereas for $\rho > 0.75$, $\delta = 2\delta$ [15]. The gain coefficient $\rho$ is also used to accept or reject the candidate designs obtained from (3). The algorithm is terminated when $\delta^{(i+1)} < \varepsilon$, or $\|x^{(i+1)} - x^{(i)}\| < \varepsilon$ (here, $\varepsilon = 10^{-3}$). For more detailed discussion on TR-based optimization of problems represented in terms of response features, see [12], [13].

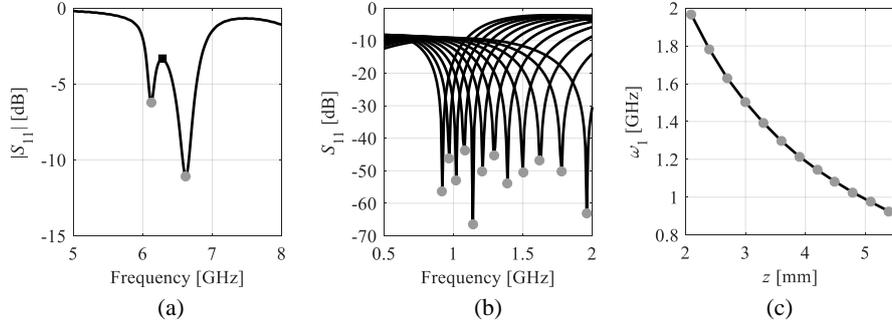

**Fig. 2.** Extraction of the response features: (a) level $S_1$ (■) at the local maximum and frequencies $\omega_1$ and $\omega_2$ (●) at a local minimum, (b) family of frequency responses obtained as a function of $z$, and (c) feature based response that corresponds to resonant frequencies along $z$. Note less non-linear changes of feature coordinates compared to frequency characteristics.

## 4    Numerical Results

The proposed topology-agnostic antenna is optimized using the method of Section 3. The initial design $x^{(0)}$ = [30 0.13 2.62 0.22 0.36 0.48 0.57 0.59 0.56 0.53 0.44 0.34 0.36 0.35 0.43 0.52 0.54 0.38 0.29 0.43 0.42 0.42 0.47 0.5 0.57 0.41 0.29 0.22 0 0.17 0.18 0.26 0.26 0.26 0.29 0.25 0.3 0.29 0.42 0.3 0.31 0.05 0.05 0.31 0.13 0.43 0.36 0.31 0.29 0.04 0.14 0.3 0.6]$^T$ is selected from a set of 200 randomly generated candidate solutions evaluated using the min-max classifier (2). The objective function for TR-based optimization (note that FD perturbations are set as 2% of $x^{(i)}$) is given as:

$$U(x) = \sum \left( \max\left( S_{\{1,2\}}(x) - S_t, 0 \right)^2 \right) + \beta \left\| \begin{array}{c} \omega_3(x) - \omega_{t.1} \\ \omega_4(x) - \omega_{t.2} \end{array} \right\| \tag{5}$$

The function (5) is calculated using four feature points: $S_1$ and $S_2$ refer to the response levels at local maxima within the bandwidth of interest, whereas $\omega_3$ and $\omega_4$ denote the local minima of the response that are close to the target frequencies $\omega_{t.1}$ = 6.2 GHz and $\omega_{t.2}$ = 6.8 GHz. The threshold on the reflection is $S_t$ = –10.2 dB, whereas $\beta$ = 100 (cf. Fig. 2). The optimized design $x^*$ = [30.04 0.13 2.47 0.21 0.36 0.47 0.56 0.61 0.58 0.52 0.46 0.34 0.36 0.31 0.44 0.53 0.55 0.38 0.33 0.44 0.41 0.42 0.49 0.52 0.57 0.41 0.28 0.23 0.01 0.16 0.17 0.25 0.26 0.27 0.29 0.24 0.3 0.29 0.42 0.29 0.29 0.08 0.07 0.32 0.13 0.43 0.36 0.31 0.29 0.04 0.15 0.32 0.61]$^T$ is found after 10 TR iterations (total cost: 221 EM simulations; see Fig. 1 for antenna visualization). It features –10 dB reflection within the operational range from 6.12 GHz to 6.84 GHz. The response at the initial and optimized designs, as well as the radiation pattern (in xz-plane) at 6.5 GHz—center frequency for 5th channel of UWB spectrum—are shown in Fig. 3. The antenna is characterized by a dual-lobe far-field response with local minimum and maximum at 0° and 35° angles, respectively. The realized gain for the latter is around 7 dB. It is worth noting that the obtained unorthodox radiation pattern is a by-product of the reflection-oriented optimization. The characteristics make the antenna suitable for installation in convex corners within the facilities dedicated for in-door monitoring using UWB-based systems dedicated [8].

The antenna of Section 2 has also been optimized using a standard TR-based algorithm with min-max design objective given as $U(x) = \max(R(x))_{\omega t.1 \leq f \leq \omega t.2}$. The benchmark method has been terminated after 7 iterations due to lack of $U(x)$ improvement with the final solution being the same as an initial design $x^{(0)}$. The obtained results suggest that, for the considered design problem, feature-based optimization is an useful tool for determination of satisfactory antenna geometries.

## 5    Conclusion

In this work, a bi-stage specification-oriented design of topology-agnostic planar antenna has been considered. The method involves identification of the starting point from a set of randomly generated designs using a classifier based on a simple min-max metric followed by its local tuning in a feature-based trust-region framework. The 53-dimensional topology has been generated at a cost of 421 EM simulations. The final design is characterized by –10 dB reflection within 6.2 GHz to 6.8 GHz band and a dual-lobe radiation pattern with minima in the direction perpendicular to the radiator and a maxima at the 35° angle. The structure is applicable for in-door positioning in challenging propagation conditions. The considered feature-enhanced TR design has been favorably benchmarked against a standard TR optimization performed in the frequency-domain. Future work will focus on implementation of the radiator as a component of the localization system.

**Acknowledgments.** This work was supported in part by the National Science Centre of Poland Grants 2020/37/B/ST7/01448 and 2021/43/B/ST7/01856, National Centre for Research and Development Grant NOR/POLNOR/HAPADS/0049/2019-00, and Gdansk University of Technology (Excellence Initiative - Research University) Grant 16/2023/IDUB/IV.2/EUROPIUM.

**Disclosure of Interests.** The authors declare no conflict of interests.

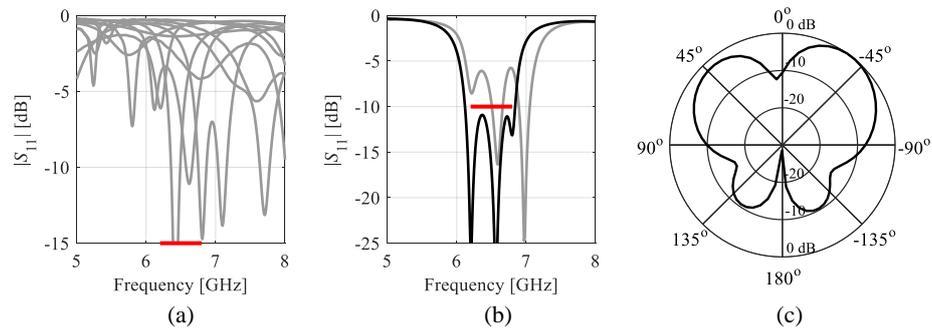

**Fig. 3.** Optimization of topology-agnostic antenna: (a) a few of randomly selected designs and bandwidth of interest for the min-max classifier (red line), (b) frequency responses at $x^{(0)}$ and $x^*$, as well as (c) xz-plane (cf. Fig. 1) radiation pattern at 6.5 GHz.